\documentclass[12pt]{article}
\usepackage{framed,epsf,latexsym,amsmath,
amssymb,amscd,mathrsfs}
\usepackage{enumitem}                                 % The preamble begins here.
\usepackage{graphicx}
\usepackage{algorithm}
\usepackage{latexsym}
\usepackage{algorithmic}
\usepackage[pdftex,colorlinks]{hyperref}
\usepackage[margin=1.0in]{geometry}
\usepackage{color}
\usepackage{amsmath,amssymb,graphicx}

\usepackage[round,comma]{natbib}
\newtheorem{theorem}{Theorem}[section]

\newtheorem{Counter-example}[theorem]{Counter example}

\newtheorem{Corollary}[theorem]{Corollary}

\newtheorem{Conjecture}[theorem]{Conjecture}

\newcommand{\ignore}[1]{}

\newcommand\blfootnote[1]{%
  \begingroup
  \renewcommand\thefootnote{}\footnote{#1}%
  \addtocounter{footnote}{-1}%
  \endgroup
}

\title{Affine Embeddings of Cantor Sets on the Line}

\author{Amir Algom}

\begin{document}
\maketitle

\begin{abstract}
Let\blfootnote{Supported by ERC grant 306494}\blfootnote{\emph{2010 Mathematics Subject Classication}. 28A80, 37C45, 11B30} $s\in (0,1)$, and let $F\subset \mathbb{R}$ be a self similar set such that $0 < \dim_H F \leq s$ . We prove that there exists $\delta= \delta(s) >0$ such that if $F$ admits an affine embedding into a homogeneous self similar set $E$ and $0 \leq \dim_H E - \dim_H F < \delta$ then (under some mild conditions on $E$ and $F$) the contraction ratios of $E$ and $F$ are logarithmically commensurable. This provides more evidence for Conjecture 1.2 of \cite{feng2014affine}, that states that these contraction ratios are logarithmically commensurable  whenever $F$ admits an affine embedding into $E$ (under some mild conditions). Our method is a combination of an argument based on the approach of \cite{feng2014affine} with a new result by \cite{hochman2016some}, which is related to the increase of entropy of measures under convolutions.
\end{abstract}

\section{Introduction}
Let $A,B \subset \mathbb{R}$. We say that $A$ can be affinely embedded into $B$ if there exists an affine map $g:\mathbb{R}\rightarrow \mathbb{R}, \quad g(x)=\gamma \cdot x +t$ for $\gamma\neq 0, t\in \mathbb{R}$, such that $g(A)\subseteq B$. The motivating problem of this paper to show that if $A$ and $B$ are two Cantor sets and $A$ can be affinely embedded into $B$, then the there should be some arithmetic dependence between the scales appearing in the  multiscale  structures of $A$ and $B$. This phenomenon is known to occur e.g. when $A$ and $B$ are two self similar sets satisfying the strong separation condition that are Lipschitz equivalent, see \cite{falconer1992lipschitz}. Similarly, one expects that if $g$ is an affine self embedding of $A$ into $A$ (i.e. $g(A) \subseteq A$) then the similarity ratio of $g$ should be  arithmetically dependent on the scales appearing in the  multiscale  structure of $A$. Evidence of the latter phenomenon is given e.g. by  Logarithmic Commensurability Theorem (Theorem 1.1 from \cite{feng2009structures}), and the results of \cite{hochman2010geometric}.

 More specifically, we are interested in studying the following Conjecture, formulated by \cite{feng2014affine}. We remark that the exact definitions of the notions we discuss are given in Section 2.  Let $F$ and $E$ be two self similar sets in $\mathbb{R}^d$, generated by IFS's $\Phi = \lbrace \phi_i \rbrace_{i=1} ^l$ and $\Psi = \lbrace \psi_j  \rbrace_{j=1} ^m$, respectively. Denote, for $1\leq i \leq l$ and $1 \leq j \leq m$, the contraction ratio of $\phi_i$ by $\alpha_i \in (0,1)$ and of $\psi_j$ by $\beta_j \in (0,1)$. 
\begin{Conjecture} \cite{feng2014affine} \label{Fengs conjecture}
Suppose that $E$ and $F$ are totally disconnected and that $F$ can be affinely embedded into $E$. Then for every $1 \leq i \leq l$ there exists $t_{i,j} \in \mathbb{Q}, t_{i,j} \geq 0$ such that 
\begin{equation*}
\alpha_i = \Pi_{j=1} ^m \beta_j ^{t_{i,j}}
\end{equation*}
In particular, if $\beta_j = \beta$ for all $1\leq j \leq m$ then for every $1 \leq i \leq l$
\begin{equation*}
\frac{\log \alpha_i}{\log \beta} \in \mathbb{Q}.
\end{equation*}
\end{Conjecture}
Note that these arithmetic conditions on $\alpha_i,\beta_j$ hold true when $\Phi$ and $\Psi$ satisfy the Strong Separation condition, and $E$ and $F$ are Lipschitz equivalent, by \cite{falconer1992lipschitz}. However, for Lipschitz embeddings no arithmetic conditions are required. Indeed, \cite{deng2011bilipschitz} proved if $E$ and $F$ are attractors of IFS's $\Phi$ and $\Psi$ that satisfy the Strong Separation condition and $\dim_H F < \dim_H E$ then $F$ can be Lipschitz embdded into $E$ (where $\dim_H F$ is the Hausdorff dimension of $F$).

\cite{feng2014affine} exhibited several cases when the Conjecture holds (and we shall soon recall them). The main result proved in this paper provides more evidence  for Conjecture \ref{Fengs conjecture}. From this point and throughout the rest of the paper,  $F$ and $E$ are two self similar sets, that admit generating IFS's $\Phi$ and $\Psi$, respectively.

\begin{theorem} \label{Main Theorem}
For every $0<s<1$ there exists $\delta=\delta (s) >0$ such that:
if $0 < \dim_H F \leq s$, $\Phi$ satisfies the OSC, $\Psi$ is homogeneous and  satisfies the SSC, and 
\begin{equation*}
0 \leq \dim_H E - \dim_H F < \delta
\end{equation*}
then $F$ and $E$ satisfy the statement of Conjecture \ref{Fengs conjecture}.
\end{theorem}
We remark that we have not worked out an effective estimate on the size of $\delta$, though this is possible.

As we mentioned earlier, \cite{feng2014affine} proved several cases when Conjecture \ref{Fengs conjecture} holds, which we now wish to discuss. We denote by $C_\rho$ the attractor of the IFS $\lbrace \rho\cdot x, \rho\cdot x + 1-\rho \rbrace$ (this is the central $\rho$ Cantor set). Recall that a Pisot number is an algebraic integer $>1$ whose algebraic conjugates are all inside the unit disc. For example, all integers greater than $1$ are Pisot numbers, as is $\sqrt{2}-1$. The following was proved by \cite{feng2014affine}:
\begin{enumerate}
\item Conjecture \ref{Fengs conjecture} is valid when $E$ has a generating IFS $\Psi$ with the strong separation condition that has uniform contraction ratio $\beta$, and $\dim_H E < \frac{1}{2}$ (no conditions on $F$, or on the contraction ratios $\alpha_i, 1\leq i \leq l$ associated with its generating IFS $\Phi$, are required).

\item  If $\alpha\leq  \beta < \sqrt{2}-1$ then Conjecture \ref{Fengs conjecture} is valid for $C_\alpha$ and $C_\beta$.

\item Assume $\theta = \frac{1}{\alpha}$ is a Pisot number $>2$. Let $E \subset \mathbb{R}$ be a self similar set generated by the IFS $\Psi = \lbrace \psi_i \rbrace_{j=1} ^m$, $2\leq m < \theta$, that has uniform contraction ratio  $\beta \in \mathbb{Z}(\theta)$ (the ring of $\theta$ over $\mathbb{Z}$). Then Conjecture \ref{Fengs conjecture} is valid for any self similar set $F$. 

Moreover, if $F$ is a self similar set generated by an IFS with contraction ratio $\alpha_i , 1\leq i \leq l$, that can be affinely embedded into $E$, then $\frac{1}{\alpha_i}$ is a Pisot number for $i=1,...,l$. 
\end{enumerate}

The main idea behind the proof of item 1 was to show that if $F$ can be affinely embedded into $E$ but $\frac{\log \alpha_i}{\log \beta} \notin \mathbb{Q}$ for some $1\leq i \leq l$, then the set $\lbrace |x-y|: x,y\in E \rbrace$ contains a non-degenerate interval. Similarly, the first step towards the proof of the item 2 involves using the self similarity structure of $C_\alpha, C_\beta$ to prove that if $\frac{\log \alpha}{\log \beta} \notin \mathbb{Q}$ and $C_\alpha$ has an affine embedding into $C_\beta$, then for every $\lambda \in (0, \frac{1-2 \alpha}{\alpha}]$ there exist $c\in C_\beta$ such that $\lambda \cdot C_\alpha + c \subseteq C_\beta$. Item 3 uses ideas from harmonic analysis.

Let us now briefly survey our methods. The proof of our main result consists of two steps. The first is to show that, assuming the conditions of our main result, if $F$  can be affinely embedded into $E$ then there exists a set $\Sigma$ consisting of affine embeddings of $F$ into $E$, such that $\dim_H \Sigma \geq 1$ (where we identify $\Sigma$ with a subset of $\mathbb{R}^2$ - see the beginning of section 2.1). To do this, we employ an analogue of the argument used by \cite{feng2014affine} to prove case 1 above.

The second step involves an application of  Theorem 4.2 from the recent paper \cite{hochman2016some}. The variant of this Theorem relevant to the present work roughly says that if  $F$ supports a "nice" measure (e.g. a self similar measure of the same dimension), and $F$ admits a big set of affine embeddings into $E$ (e.g. of Hausdorff dimension at least $1$), then there exists a measure $\theta$ supported on $E$ with entropy dimension $> \dim_H F + \delta ( s)$, where $ \delta (s) >0$ depends only on the choice of $s$ such that $\dim_H F \leq s<1$. See exact definitions and statements in section 2.2.

Let us end this introduction by  stating two Corollaries that follow from Theorem \ref{Main Theorem}. First, if one assumes that $E = F$, then slightly modifying of the proof of Theorem \ref{Main Theorem} we obtain the following version of the Logarithmic Commensurability Theorem (Theorem 1.1 from \cite{feng2009structures}). Our version is actually weaker than the Theorem by \cite{feng2009structures}, since we only treat the case when $0<\dim_H F <1$ and $F$ admits an IFS $\Phi$ with the SSC, whereas they only require that $\Phi$ satisfies the OSC, and also treat the case $\dim_H F =1$.

\begin{Corollary} \label{LCT}
Let $F \subset \mathbb{R}$ be a self similar set that admits a homogeneous generating IFS $\Phi$, with uniform contraction ratio $\alpha$, that satisfies the SSC. Assume $0<\dim_H F <1$. If $ g (x) = \lambda\cdot x +t, \quad \lambda \neq 0, t\in \mathbb{R}$ satisfies $g (F) \subseteq F$ then
\begin{equation*}
\frac{\log \alpha}{\log |\lambda|} \in \mathbb{Q}
\end{equation*}
\end{Corollary}

We shall say that $F$ can be $C^1$ embedded into $E$ if $g(F) \subseteq E$ for some $C^1$-diffeomorphism $g:\mathbb{R} \rightarrow \mathbb{R}$. We denote the set of all $C^1$-diffeomorphisms $g:\mathbb{R} \rightarrow \mathbb{R}$ by $\text{diff}^1 (\mathbb{R})$. To derive the second Corollary, let us first recall the following Theorem, that relates the notions of affine embedding and $C^1$ embedding for self similar sets. The Theorem below is stated in the form that is relevant to this present work (but not the most general form):

\begin{theorem} \cite{feng2014affine} \label{Theorem feng C1 embedding}
Suppose both $\Phi$ and $\Psi$ satisfy the OSC. Then $F$ can be $C^1$ embedded into $E$ if and only if $F$ can be affinely embedded into $E$. Furthermore, if $F$ can not be affinely embedded into $E$ then
\begin{equation*}
\sup_{ g\in \text{diff}^1 (\mathbb{R})} \dim_H g(F) \cap E < \dim_H F
\end{equation*}
\end{theorem}

Thus, combining Theorem \ref{Theorem feng C1 embedding} with Theorem \ref{Main Theorem}, yields the following Corollary:

\begin{Corollary} \label{Thereom C1 embedding}
Suppose that $F$ and $E$ satisfy the conditions of Theorem \ref{Main Theorem}. Suppose that there exists a map $\phi \in \Phi$ with contraction ratio $\alpha$ such that $\frac{\log \alpha}{\log \beta} \notin \mathbb{Q}$, where $\beta$ is the uniform contraction ratio of maps in $\Psi$. Then
\begin{equation*}
\sup_{ g\in \text{diff}^1 (\mathbb{R})} \dim_H g(F) \cap E < \dim_H F
\end{equation*}

\end{Corollary}

\textbf{Acknowledgments} This paper is part of the author's research  towards a PhD dissertation, conducted at the Hebrew University of Jerusalem. I would like to sincerely  thank my adviser, Michael Hochman, for his continuous encouragement and support, and for many useful comments and suggestions. I would also like to thank Ariel Rapaport for many helpful discussions. The author is also
grateful for the hospitality and support received from ICERM as part of the spring 2016 program
on dimension and dynamics.

\section{Preliminaries}
\subsection{Self similar sets}
Let $G$ denote the group of similarities of the real line, and $S\subset G$ the set of contracting similarities. Then we can naturally identify $G$ and $S$ with subsets of $\mathbb{R}^2$, as
\begin{equation*}
G= \lbrace g:\mathbb{R} \rightarrow \mathbb{R}: g (x) = \alpha \cdot x + t, \quad (\alpha,t)\in \mathbb{R}^2 - \lbrace 0 \rbrace \times \mathbb{R} \rbrace
\end{equation*}
\begin{equation*}
S= \lbrace g :\mathbb{R} \rightarrow \mathbb{R}: g (x) = \alpha \cdot x + t, \quad (\alpha,t)\in (0,1)\times \mathbb{R} \rbrace.
\end{equation*}
As we are working in $\mathbb{R}$, $G$ in fact consists of all the invertible affine maps $\mathbb{R} \rightarrow \mathbb{R}$.

Let $\Phi = \lbrace \phi_i \rbrace_{k=1} ^l , l\in \mathbb{N}, l\geq 2$ be a family of contractions $\phi_i : \mathbb{R}^d \rightarrow \mathbb{R}^d, d\geq 1$. The family $\Phi$ is called an iterated function system, abbreviated IFS, the term being coined by \cite{hutchinson1981fractals}, who defined them and studied some of their fundamental properties. In particular, he proved that there exists a unique compact $\emptyset \neq F \subset \mathbb{R}^d$ such that $F = \bigcup_{i=1} ^l \phi_i (F)$. $F$ is called the attractor of $\Phi$, and $\Phi$ is called a generating IFS for $F$.

A set $F \subset \mathbb{R}$ will be called self similar if there exists a generating IFS $\Phi$ for $F$ such that $\Phi \subset S$.  If this IFS has the additional property that all the maps $\phi \in \Phi$ share the same contraction ratio, then we shall say that $\Phi$ is a homogeneous generating IFS for $\Phi$. We shall always assume, throughout this paper, that our self similar sets satisfy $0<\dim_H F <1$, where $\dim_H F$ denotes the Hausdorff dimension of $F$.

For an IFS $\Phi = \lbrace \phi_i \rbrace_{i=1} ^l$ and its attractor $F$, a cylinder set is a set of the form $\phi_{i_1} \circ ... \circ \phi_{i_k} (F)$, where $\phi_i \in \Phi$ for all $i$ and $k\in \mathbb{N}$. Writing $I= (i_1,...,i_k)\in \lbrace 0,...,l \rbrace^k$, we use the notation $\phi_{i_1} \circ ... \circ \phi_{i_k} = \phi_I$. Thus, cylinder sets have the form $\phi_I (F), I\in \lbrace 0,...,l \rbrace^*$, and $|I|$ denotes the length of the word $I$.

We shall say that an IFS $\Phi$ satisfies the open set condition, abbreviated OSC, if there exists some $\emptyset \neq U \subset \mathbb{R}^d$ such that $U$ is open, $i\neq j \Rightarrow \phi_i (U) \cap \phi_j (U) = \emptyset$ and $U \supset \bigcup_{i=1} ^l \phi_i (U)$. We shall say that an IFS $\Phi$ satisfies the strong separation condition, abbreviated SSC, if  $i\neq j \Rightarrow \phi_i (F) \cap \phi_j (F) = \emptyset$. It is straightforward to verify that if $\Phi$ satisfies the SSC then it also satisfies the OSC.

\subsection{Growth of entropy dimension of measures}
The objective of this subsection is to state some of the results of the upcoming paper \cite{hochman2016some}, that are key to the proof of Theorem \ref{Main Theorem}. We do not aim at explaining why these results are true; we leave that to \cite{hochman2016some}.  We only quote  what we need and proceed to the proof of Theorem \ref{Main Theorem}, using these results as a black box.

Let $P(X)$ denote the space of Borel probability measure supported on $X\subset \mathbb{R}^d$. We first recall the definition of entropy dimension of a measure.    Let 
\begin{equation*}
D_n = \lbrace [\frac{k}{2^n}, \frac{k+1}{2^n}) \rbrace_{k\in \mathbb{Z}}
\end{equation*}
denote the level $n$ dyadic partition of $\mathbb{R}$. Let 
\begin{equation*}
H(\theta, \mathcal{E}) =  -\sum_{E\in \mathcal{E}} \theta (E) \log \theta (E)
\end{equation*}
denote the Shannon entropy of a probability measure $\theta \in P(\mathbb{R})$ with respect to a partition $\mathcal{E}$ of $\mathbb{R}$. Then the entropy dimension of $\theta$ is defined as
\begin{equation*}
\dim_e \theta = \lim_{n\rightarrow \infty} \frac{1}{n} H(\theta, D_n).
\end{equation*}
If the above limit does not exist, we define the upper entropy dimension $\overline{\dim}_e \theta$ by taking $\limsup$. Note that if $\theta$ is a self similar measure then this limit exists (see e.g. \cite{hochman2016some} Proposition 2). We also note that if $\theta \in P(\mathbb{R})$ is supported on a set $Y$ then
\begin{equation} \label{Equation dimension}
\overline{\dim}_B Y \geq \overline{\dim}_e \theta \geq \underline{\dim}_H \theta 
\end{equation}
where $\overline{dim}_B Y$ is the upper box dimension of $Y$, and
\begin{equation*}
\underline{\dim}_H \theta  =  \inf \lbrace \dim A : A \text{ is Borel }, \theta (A)>0 \rbrace 
\end{equation*}

Note that if $F$ is the attrator of an IFS that satisfies the OSC then there exists a self similar measure $\mu$ supported on $F$ such that  $\underline{\dim}_H \mu = \dim_H F$. We shall call this measure a self similar measure of maximal dimension. See \cite{falconer1986geometry},\cite{mattila1999geometry}, \cite{bishop2013fractal} for more details. 

Next, let $v\in P(G)$ and $\mu \in P(\mathbb{R})$ be compactly supported probability measures. Then $\nu.\mu \in P(\mathbb{R})$ is defined as the push forward of $\nu \times \mu$ via the action map $(g,x) \mapsto g (x)$, from $G\times R$ to $\mathbb{R}$. Note that this is a smooth map defined on an open subset of $\mathbb{R}^2 \times \mathbb{R}$, so $\nu.\mu$ is a Borel probability measure on $\mathbb{R}$. This concept brings us to Theorem 6  from \cite{hochman2016some} (which is actually a consequence of the more general results from \cite{hochman2015self}), that shall be used as a black box. Note that we state the version of the Theorem that we require, but certainly not the most general form.

\begin{theorem} \cite{hochman2016some} \label{Theorem 4.2}
Let $s\in (0,1)$, then there exists some $\delta=\delta(s) >0$ such that:

Let $\mu \in P(\mathbb{R})$ and $\nu \in P(G)$ be compactly supported measures. Suppose that $\mu$ is a  self similar measure of maximal dimension with respect to some IFS that satisfies the OSC, and that its attractor $F$ satisfies $0 < \dim_H F \leq  s $. Suppose in addition that $\overline{\dim}_e \nu \geq  1$.  Then
\begin{equation*}
\overline{\dim}_e \nu. \mu \geq \dim_H F + \delta.
\end{equation*} 
\end{theorem}

We remark that to obtain Theorem \ref{Theorem 4.2} from Theorem 6 of \cite{hochman2016some}, one assumes the conditions of Theorem \ref{Theorem 4.2} and plugs into Theorem 6
\begin{equation*}
\epsilon = 1-s.
\end{equation*}

\section{Proof of Theorem \ref{Main Theorem} and Corollary \ref{LCT}}
\textbf{Proof of Theorem \ref{Main Theorem}} Let $s\in (0,1)$, and let $F$ and $E$ be real self similar sets, generated by the IFS's  $\Phi$ and $\Psi$ respectively. Assume $\Psi$  satisfies the SSC and is homogeneous, and that $\Phi$ satisfies the OSC. Denote the  contraction ratios of the maps in $\Phi$ by $\alpha_i, 1\leq i \leq l$, and the uniform contraction ratio of the maps in $\Psi$ by $\beta$. Assume $0 <\dim_H F \leq s<1$ and let $\delta = \delta (s)>0$ be as in the conclusion of  Theorem \ref{Theorem 4.2}, and assume
\begin{equation} \label{Equation dimension of E and F}
0\leq \dim_H E - \dim_H F < \delta. 
\end{equation}

Suppose that $g(F) \subseteq E$ is an affine embedding. Suppose, towards a contraction,  that we have $\frac{\log \alpha_i }{\log \beta} \notin \mathbb{Q}$, for some $1\leq i \leq l$. Without the loss of generality, suppose $i=1$. To make notation easier, denote $\alpha =  \alpha_1$, so that  $\frac{\log \alpha }{\log \beta} \notin \mathbb{Q}$.

\textbf{Step 1} We first claim that there exists a compact set $\Sigma \subset G$ of affine embeddings of $F$ into $E$ such that $p_1 (\Sigma)$ contains an interval, where $p_1 :\mathbb{R}^2 \rightarrow \mathbb{R}$ is the projection $p_1 (x,y)=x$ (recall that we are identifying $G$ with a subset of $\mathbb{R}^2$). Write $g(x)= \gamma \cdot x +b$. Let 
\begin{equation*}
\kappa = \min_{i\neq j} d( \psi_i (E), \psi_j (E))
\end{equation*}
then $\kappa >0$ by the assumption that $\Psi$ has the SSC. Denote  $c = \gamma \cdot \text{diam} (F)$. Fix $p,N$ such that $\beta ^p < \frac{\kappa}{c}$ and $\frac{\log \alpha}{\log \beta} \cdot N >p$.

Then for all $n \in \mathbb{N}$ with $n>N$ we have 
\begin{equation*}
g( \phi_{1^n} (F)) \subset g(F) \subseteq E
\end{equation*}
Note that 
\begin{equation*}
\phi_{1^n} (F) = \alpha^n \cdot F + e_n , \quad e_n \in \mathbb{R}.
\end{equation*}
Hence we have 
\begin{equation*}
g( \phi_{1^n} (F)) =  \gamma \cdot \alpha^n  \cdot F + \gamma \cdot e_n +b \subset E. 
\end{equation*}
Let $l_n = [ \frac{\log \alpha}{\log \beta} \cdot n]$, where for $x\in \mathbb{R}$, $[x]$ denotes its floor integer value. Then 
\begin{equation} \label{equation diameter}
\text{diam} ( \gamma \cdot \alpha^n  \cdot F + \gamma \cdot e_n +b ) = \alpha^n \cdot \gamma \cdot \text{diam}(F) = \alpha^n \cdot c = \beta^{\frac{\log \alpha}{\log \beta} \cdot n} \cdot c 
\end{equation}
\begin{equation*}
=  \beta^{\frac{\log \alpha}{\log \beta} \cdot n - p} \cdot \beta^p \cdot c < \beta^{\frac{\log \alpha}{\log \beta} \cdot n - p} \cdot \kappa \leq \beta^{l_n -p } \cdot \kappa < \beta^{l_n -p-1 } \cdot \kappa.
\end{equation*}

Note that we have
\begin{equation*}
E = \bigcup_{I\in \lbrace 1,...,m \rbrace^{l_n - p}} \psi_I (E),
\end{equation*}
and by the definition of $\kappa$,
\begin{equation*}
 \min_{I\neq J, I,J \in \lbrace 1,...,m \rbrace^{l_n - p}} d( \psi_I (E), \psi_J (E)) = \beta^{l_n - p-1} \cdot  \kappa .
\end{equation*}
Therefore,  by equation \eqref{equation diameter} and since 
\begin{equation*}
 \gamma \cdot \alpha^n  \cdot F + \gamma \cdot e_n +b \subset E
\end{equation*}
it follows that $ \gamma \cdot \alpha^n  \cdot F + \gamma \cdot e_n +b$ intersects  a unique $\psi_I (E)$ for $|I| = l_n -p, I\in \lbrace 1,...,m \rbrace^{l_n - p}$. Therefore, for some $r_n \in \mathbb{R}$
\begin{equation*}
 \gamma \cdot \alpha^n  \cdot F + \gamma \cdot e_n +b \subset \psi_I (E) = \beta^{l_n - p} \cdot E +r_n.
\end{equation*}
Finally, we obtain
\begin{equation} \label{Equation embedding}
\beta^{p+ \lbrace \frac{\log \alpha}{\log \beta} \cdot n \rbrace } \cdot \gamma \cdot F +\frac{\gamma \cdot e_n +b -r_n}{\beta^{l_n -p}} \subset E
\end{equation}
where for $x\in \mathbb{R}$ we define $\lbrace x \rbrace$ as the fractional part of $x$ (so that $x=[x]+\lbrace x \rbrace$).

Taking limits with respect to $n$ in  equation \eqref{Equation embedding}, and recalling our assumption towards a contradiction that $\frac{\log \alpha}{\log \beta} \notin \mathbb{Q}$, we claim that for every $\eta \in [\beta^{p+1} \cdot \gamma, \beta^{p} \cdot \gamma]$ there exists some $t \in \mathbb{R}$ such that $\eta \cdot F + t \subseteq E$ (where we assume without the loss of generality that $\gamma >0$, otherwise $\eta \in [\beta^{p} \cdot \gamma, \beta^{p+1} \cdot \gamma]$ and we proceed in a similar manner). Moreover, the set of these $t$'s that we obtain is bounded (as a subset of $\mathbb{R}$).

This can be seen by fixing some element $f\in F$, and observing equation \eqref{Equation embedding} along a sequence $n_k$ such that $ \beta^{p+ \lbrace \frac{\log \alpha}{\log \beta} \cdot n_k \rbrace } \rightarrow \eta$. We obtain
\begin{equation*}
\frac{\gamma \cdot e_{n_k} +b -r_{n_k}}{\alpha^{l_{n_k} -p}} \in E - \beta^{p+ \lbrace \frac{\log \alpha}{\log \beta} \cdot n_k \rbrace } \cdot f
\end{equation*}
It follows that $\frac{\gamma \cdot e_{n_k} +b -r_{n_k}}{\alpha^{l_{n_k} -p}}$ is a  bounded sequence as $k\rightarrow \infty$, and hence has a converging sub-sequence, with limit that shall be denoted $t$. Taking the limit in \eqref{Equation embedding} along this sub-sequence, we see that $\eta \cdot F + t \subseteq E$. That the set  $t$'s that we obtain is bounded follows since by this proof, these are elements in the interval
\begin{equation*}
[\min E - \beta^p \cdot \gamma \cdot \max F, \max E - \beta^{p+1} \cdot \gamma \cdot \min F].
\end{equation*}

Let $A \subset G$ denote the set of the embeddings that we obtain by this procedure. Then $A$ is a bounded subset of $G$, which is a subset of $\mathbb{R}^2$, and the projection of $A$ onto the $x$-axis (which corresponds to the similarity ratios) is the interval $[\beta^{p+1} \cdot \gamma, \beta^{p} \cdot \gamma]$. Therefore, $\Sigma = cl (A)$ (the closure of $A$) is a compact subset of $G$ such that $p_1 (\Sigma) \supseteq p_1 (A)$,  so $p_1 (\Sigma)$ contains an interval. This finishes the proof of step 1.

\textbf{Step 2} We now apply the machinery of \cite{hochman2016some}. Let $\mu$ be a self similar measure of maximal dimension supported on $F$ (recall that we are assuming that $\Phi$, the generating IFS for $F$, satisfies the OSC). Let $\theta$ be the (normalized) Lebesgue measure on the interval we found to be contained in $p_1 (\Sigma)$. Apply Theorem 1.20 from \cite{mattila1999geometry} to  obtain a compactly supported measure $\nu \in P(\Sigma)$  such that $p_1 \nu = \nu \circ p_1 ^{-1} = \theta$. Therefore,  $\underline{\dim}_H \nu \geq  1$. Thus,  By equation \eqref{Equation dimension}, $\overline{\dim}_e \nu \geq  1$.  Then $\nu.\mu$ is a probability measure supported on the image of $\Sigma \times F$ under the action map $(g,x) \mapsto g(x)$. Therefore, $\nu.\mu$ is supported on $E$. Since $\mu$ is a self similar measure of maximal dimension supported on $F$, $0<\dim_H F \leq s$, and $\overline{\dim}_e \nu \geq  1$, we may apply Theorem \ref{Theorem 4.2}. We thus obtain that 
\begin{equation*}
\overline{\dim}_B E \geq \overline{\dim}_e \nu.\mu  \geq \dim_H F + \delta,
\end{equation*}
where $\delta = \delta (s)>0$, by another application of equation \eqref{Equation dimension}. 

However, by e.g. \cite{feng2009dimension} Theorem 2.13, we have that since $E$ is a self similar set,
\begin{equation*}
\dim_H E = \dim_B E = \overline{\dim}_B E \geq \dim_H F + \delta
\end{equation*}
This implies that $\dim_H E - \dim_H F \geq \delta$, contradicting equation \eqref{Equation dimension of E and F}. This proves that $\frac{\log \alpha_i }{\log \beta} \in \mathbb{Q}$ for all $i$, concluding the proof. \hfill{$\Box$}

$$ $$

\textbf{Proof of Corollary \ref{LCT}} Assume now that $E=F$ is a self similar set generated by the homogeneous IFS $\Phi$ with the SSC and uniform contraction ratio $\alpha$. Suppose $g(x)=\gamma \cdot x +b$ is an affine map $g:\mathbb{R} \rightarrow \mathbb{R}$ such that $g(F) \subseteq F$. We prove that  $\frac{\log \alpha}{\log |\gamma|} \in \mathbb{Q}$.

Suppose towards a contradiction that $\frac{\log \alpha}{\log |\gamma|} \notin \mathbb{Q}$.  The idea is to follow in the lines of the proof of Theorem \ref{Main Theorem}, showing that (as in step 1) this implies that we can find a compact subset $\Sigma \subset G$ such that $p_1 (\Sigma)$ contains an interval, and for every $g\in \Sigma$, $g(F) \subseteq F$. Then we apply  (as in step 2) Theorem \ref{Theorem 4.2} to obtain that $\overline{\dim}_B F \geq \dim_H F +\delta$, where $\delta = \delta (\dim_H F)>0$. Since $F$ is self similar, $\dim_H F = \overline{\dim}_B F$, thus we have $ \dim_H F \geq \dim_H F +\delta$, a contradiction.

For brevity, we sketch the proof of the fact that we can find a compact subset $\Sigma \subset G$ such that $p_1 (\Sigma)$ contains an interval and for every $g\in \Sigma$, $g(F) \subseteq F$. We may assume $\gamma>0$, since otherwise the similarity ratio of $g^2 = g\circ g$ is $\gamma^2>0$, and $g^2 (F) \subseteq g(F) \subseteq F$. It is clearly sufficient to prove the result for $\gamma^2>0$, so we assume without the loss of generality that $\gamma>0$.

Let 
\begin{equation*}
\kappa = \min_{i\neq j} d( \phi_i (F), \phi_j (F)),
\end{equation*}
so $\kappa >0$ by the assumption that $\Phi$ has the SSC. Denote  $c = \gamma \cdot \text{diam} (F)$. Fix $p,N$ such that $\alpha ^p < \frac{\kappa}{c}$ (note that we may assume $|\gamma| <1$, since $|\gamma|=1$ already contradicts $\frac{\log \alpha}{\log |\gamma|} \notin \mathbb{Q}$) and $\frac{\log \gamma}{\log \alpha} \cdot N >p$.

Then for all $n \in \mathbb{N}$ with $n>N$ we have, for some $b_{n+1} \in \mathbb{R}$ where $b_1 = b$,
\begin{equation*}
g^{n+1} ( F) = \gamma^{n+1} F+b_{n+1}  \subseteq F.
\end{equation*}
Let $l_n = [ \frac{\log \gamma}{\log \alpha} \cdot n]$. Then 
\begin{equation*}
\text{diam} ( \gamma^{n+1} \cdot F +b ) = \gamma^{n} \cdot \gamma \cdot \text{diam}(F)= \gamma^n \cdot c <  \alpha^{l_n - p} \kappa < \alpha^{l_n - p-1} \kappa.
\end{equation*}
By the definition of $\kappa$, we see that $ \gamma^{n+1}  \cdot F  +b$ intersects  a unique $\phi_I (F)$ for $|I| = l_n -p, I\in \lbrace 1,...,l \rbrace^{l_n - p}$. Therefore, for some $r_n \in \mathbb{R}$
\begin{equation*}
 \gamma^{n+1}  \cdot F  +b_{n+1} \subset \phi_I (F) = \alpha^{l_n - p} \cdot F +r_n.
\end{equation*}
Finally, we obtain
\begin{equation} \label{Equation embedding2}
\gamma \cdot \alpha^{p+ \lbrace \frac{\log \gamma}{\log \alpha} \cdot n \rbrace }  \cdot F +\frac{ b_{n+1} -r_n}{\alpha^{l_n -p}} \subset F.
\end{equation}

It follows that, by taking limits in \eqref{Equation embedding2} and recalling the assumption $\frac{\log \alpha}{\log \gamma} \notin \mathbb{Q}$, that there exists a subset $\Sigma \subset G$ such that for every $g\in \Sigma$ we have $g(F) \subseteq F$, and the projection of $\Sigma$ to the $x$-axis is the interval $[\gamma \cdot \alpha^{p+1}, \gamma\cdot\alpha^p]$. 

As explained in the beginning of this proof, applying Theorem \ref{Theorem 4.2}, we obtain $\dim_H F \geq \dim_H F +\delta$ for $\delta = \delta (\dim_H F)>0$, a contradiction. This concludes the proof. \hfill{$\Box$}

\section{Concluding remarks}

\begin{enumerate}
\item It is natural to ask what can one say about conjecture \ref{Fengs conjecture} in higher dimension. Similarly, one may ask (see Open Question 2 in \cite{feng2009structures}) about a generalization of Corollary \ref{LCT} in higher dimension for self affine sets (for self similar sets there is a corresponding result in \cite{elekes2010self} in $\mathbb{R}^d$).  This is a more delicate situation, though it is tractable in some cases. We shall address these problems in future works.

\item Another natural question that arises by observing Conjecture \ref{Fengs conjecture} is the following: Let $E$ and $F$ be totally disconnected self similar sets. Can one give a sufficient condition that ensures that there exists an affine embedding of $F$ into $E$?
\end{enumerate}

\bibliography{bib}  
\end{document}